\documentstyle[12pt]{article}

\font\twlgot =eufm10 scaled \magstep1 \font\egtgot =eufm8
\font\sevgot =eufm7 \font\twlmsb =msbm10 scaled \magstep1
\font\egtmsb =msbm8 \font\sevmsb =msbm7

\newfam\gotfam
\def\pgot{\fam\gotfam\twlgot}
\textfont\gotfam\twlgot \scriptfont\gotfam\egtgot
\scriptscriptfont\gotfam\sevgot
\def\got{\protect\pgot}
\newfam\msbfam
\textfont\msbfam\twlmsb \scriptfont\msbfam\egtmsb
\scriptscriptfont\msbfam\sevmsb
\def\Bbb{\protect\pBbb}
\def\pBbb{\relax\ifmmode\expandafter\Bb\else\typeout{You cann't use
Bbb in text mode}\fi}
\def\Bb #1{{\fam\msbfam\relax#1}}

\newcommand{\gd}{{\got d}}

\newcommand{\ccG}{{\got g}}

\def\thebibliography#1{\section*{References}\list
  {[\arabic{enumi}]}{\settowidth\labelwidth{#1}\leftmargin\labelwidth
    \advance\leftmargin\labelsep
    \usecounter{enumi}}
    \def\newblock{\hskip .11em plus .33em minus .07em}
    \sloppy\clubpenalty4000\widowpenalty4000
    \sfcode`\.=1000\relax}

\def\op#1{\mathop{\fam0 #1}\limits}

\newcommand{\Ker}{{\rm Ker\,}}

\newcommand{\nm}[1]{|{#1}|}
\newcommand{\beq}{\begin{equation}}
\newcommand{\eeq}{\end{equation}}
\newcommand{\ben}{\begin{eqnarray}}
\newcommand{\een}{\end{eqnarray}}
\newcommand{\be}{\begin{eqnarray*}}
\newcommand{\ee}{\end{eqnarray*}}
\newcommand{\bea}{\begin{eqalph}}
\newcommand{\eea}{\end{eqalph}}

\newcommand{\cL}{{\cal L}}

\newcommand{\cE}{{\cal E}}

\newcommand{\cF}{{\cal F}}

\newcommand{\cO}{{\cal O}}

\newcommand{\bL}{{\bf L}}

\newcommand{\al}{\alpha}

\newcommand{\bt}{\beta}
\newcommand{\dl}{\delta}
\newcommand{\la}{\lambda}
\newcommand{\La}{\Lambda}
\newcommand{\f}{\phi}
\newcommand{\om}{\omega}

\newcommand{\m}{\mu}

\newcommand{\g}{\gamma}
\newcommand{\G}{\Gamma}
\newcommand{\th}{\theta}

\newcommand{\vt}{\vartheta}

\newcommand{\up}{\upsilon}

\newcommand{\di}{{\rm dim\,}}

\newcommand{\si}{\sigma}
\newcommand{\Si}{\Sigma}
\newcommand{\w}{\wedge}

\newcommand{\ol}{\overline}
\newcommand{\dr}{\partial}
\newcommand{\ar}{\op\longrightarrow}
\newcommand{\ot}{\otimes}
\newcommand{\ap}{\approx}
\newcommand{\ve}{\varepsilon}
\newcommand{\e}{\epsilon}

\let\ssection=\section
\renewcommand{\section}{\setcounter{equation}{0}\ssection}

\newcounter{eqalph}[section]
\newcounter{equationa}[section]
\newcounter{example}[section]
\newcounter{remark}[section]
\newcounter{theorem}[section]
\newcounter{proposition}[section]
\newcounter{lemma}[section]
\newcounter{corollary}[section]
\newcounter{definition}[section]

\def\theremark{\arabic{section}.\arabic{remark}}

\def\thedefinition{\arabic{section}.\arabic{definition}}

\newenvironment{proof}{\noindent
{\bf Proof.} }{$\Box$ \medskip}
\newenvironment{rem}{\refstepcounter{remark}\medskip\noindent{\bf
Remark \theremark.} }{\medskip}
\newenvironment{ex}{\refstepcounter{remark}\medskip\noindent{\bf
Example \theremark.} }{\medskip}
\newenvironment{theo}{\refstepcounter{definition}
\bigskip\noindent{\bf Theorem \thedefinition.} \it}{\medskip}
\newenvironment{prop}{\refstepcounter{definition}
\bigskip\noindent{\bf Proposition \thedefinition.}\it}{\medskip}

\newenvironment{defi}{\refstepcounter{definition}
\bigskip\noindent{\bf Definition \thedefinition.}\it}{\medskip}

\newenvironment{eqalph}{\stepcounter{equation}
\setcounter{equationa}{\value{equation}} \setcounter{equation}{0}

\begin{eqnarray}}{\end{eqnarray}
\setcounter{equation}{\value{equationa}}}

\newcommand{\mar}[1]{}

\begin{document}
\hbox{}

{\parindent=0pt

{\large \bf Noether's second theorem in a general setting.
Reducible gauge theories}
\bigskip

{\bf D.Bashkirov}$^1$, {\bf G.Giachetta}$^2$, {\bf
L.Mangiarotti}$^2$, {\bf G.Sardanashvily}$^1$
\bigskip

\begin{small}

$^1$ Department of Theoretical Physics, Moscow State University,
117234 Moscow, Russia
\medskip

$^2$ Department of Mathematics and Informatics, University of
Camerino, 62032 Camerino (MC), Italy

\bigskip

{\bf Abstract}

We prove Noether's direct and inverse second theorems for
Lagrangian systems on  fiber bundles in the case of gauge
symmetries depending on derivatives of dynamic variables and
parameters of an arbitrary order. The appropriate notions of a
reducible gauge symmetry and Noether identity are formulated, and
their equivalence by means a certain intertwining operator is
proved.

\end{small}
\bigskip


 }

\section{Introduction}

Different variants of Noether's second theorem state that, if a
Lagrangian admits symmetries depending on parameters, its
variational derivatives obey certain relations, called the Noether
identity. In a rather general setting, this theorem has been
formulated in \cite{fulp}. Gauge symmetries and Noether identities
need not be independent, and one speaks on $N$-stage reducible
gauge symmetries and Noether identities. The notion of a reducible
Noether identity has come from that of a reducible constraint
\cite{fish}, but it involves differential operators. Note that the
conventional Batalin--Vilkovisky (BV) quantization of a classical
gauge system necessarily starts with studying an hierarchy of its
gauge symmetries and Noether identities in order to define the
multiplet of ghosts and antifields, and to construct the so-called
gauge-fixed Lagrangian \cite{barn,gom}. It should be also
emphasized that, if a gauge symmetry is reducible, the components
of a Noether current in classical field theory and the Ward
identities in quantum field theory fail to be independent.

We present Noether's second theorem and its inverse (Theorem
\ref{0510}) for Lagrangian systems on a fiber bundle $Y\to X$ in
the case of gauge symmetries depending on derivatives of dynamic
variables and parameters of an arbitrary order. Bearing in mind
the extension of the BV quantization scheme to an arbitrary base
manifold \cite{bash,cmp}, we pay the particular attention to
global aspects of Noether's second theorem. For this purpose, we
consider Lagrangian formalism on the composite fiber bundle $E\to
Y\to X$, where $E\to Y$ is a vector bundle of gauge parameters.
Accordingly, a gauge symmetry is represented by a linear
differential operator $\up$ on $E$ taking its values in the
vertical tangent bundle $VY$ of $Y\to X$.

The Noether identity for a Lagrangian $L$ is defined by a
differential operator $\Delta$ on the fiber bundle (\ref{0548})
which takes its values in the density-dual
\mar{0630}\beq
\ol E^*=E^*\op\ot_Y\op\w^n T^*X \label{0630}
\eeq
 of
$E\to Y$ and whose kernel contains the image of the
Euler--Lagrange operator $\dl L$ of $L$, i.e.,
\mar{0549}\beq
\Delta\circ \dl L=0 \label{0549}.
\eeq
Expressed in these terms, Noether's second theorem and its inverse
follow at once from the first variational formula (Proposition
\ref{g75}) and the properties of differential operators on dual
fiber bundles (Theorem \ref{0610}). Namely, there exists the
intertwining operator $\eta(\up)=\Delta$ (\ref{0616}),
$\eta(\Delta)=\up$ (\ref{0617}) such that
\mar{0634,22}\ben
&& \eta(\eta(\up))=\up, \qquad \eta(\eta(\Delta))=\Delta,
\label{0634}\\
&& \eta(\up\circ\up')=\eta(\up')\circ\eta(\up), \qquad
\eta(\Delta'\circ\Delta)=\eta(\Delta)\circ\eta(\Delta').
\label{0622}
\een
The appropriate notions of a reducible Noether identity and gauge
symmetry  are formulated, and their equivalence with respect to
the intertwining operator $\eta$ is proved (Theorem \ref{0580}).

The following two examples aim to illustrate our exposition: (i)
the gauge theory of principal connections for which gauge transformations
need not be vertical, e.g., the topological gauge theory with the
global Chern--Simons Lagrangian and the Yang--Mills gauge theory
with a dynamic metric field, (ii) a gauge system of skew symmetric
tensor fields with a reducible gauge symmetry, e.g., the
topological BF theory.

\section{Lagrangian formalism on fiber bundles}

Lagrangian formalism on a fiber bundle $Y\to X$ is phrased in
terms of the following graded differential algebra (henceforth
GDA) \cite{ander,jmp,cmp,tak2}.

The finite order jet manifolds of $Y\to X$ form an inverse
system
\mar{5.10}\beq
X\op\longleftarrow^\pi Y\op\longleftarrow^{\pi^1_0} J^1Y
\longleftarrow \cdots J^{r-1}Y \op\longleftarrow^{\pi^r_{r-1}}
J^rY\longleftarrow\cdots. \label{5.10}
\eeq
In the sequel, the index $r=0$ stands for $Y$. Accordingly, we
have the direct system
\mar{5.7}\beq
\cO^*X\op\longrightarrow^{\pi^*} \cO^*Y
\op\longrightarrow^{\pi^1_0{}^*} \cO_1^*Y \ar\cdots \cO^*_{r-1}Y
\op\longrightarrow^{\pi^r_{r-1}{}^*}
 \cO_r^*Y \longrightarrow\cdots \label{5.7}
\eeq
of GDAs $\cO_r^*Y$ of exterior forms on jet manifolds $J^rY$ with
respect to the pull-back monomorphisms $\pi^r_{r-1}{}^*$. Its
direct limit
 $\cO_\infty^*[Y]$ is a GDA
consisting of all exterior forms on finite order jet manifolds
modulo the pull-back identification.

The projective limit $(J^\infty Y, \pi^\infty_r:J^\infty Y\to
J^rY)$ of the inverse system (\ref{5.10}) is a Fr\'echet manifold.
A bundle atlas $\{(U_Y;x^\la,y^i)\}$ of $Y\to X$ yields the
coordinate atlas
\mar{jet1}\beq
\{((\pi^\infty_0)^{-1}(U_Y); x^\la, y^i_\La)\}, \qquad
{y'}^i_{\la+\La}=\frac{\dr x^\m}{\dr x'^\la}d_\m y'^i_\La, \qquad
0\leq|\La|, \label{jet1}
\eeq
of $J^\infty Y$, where $\La=(\la_k...\la_1)$ is a symmetric
multi-index, $\la+\La=(\la\la_k...\la_1)$,  and
\mar{5.177}\beq
d_\la = \dr_\la + \op\sum_{0\leq|\La|} y^i_{\la+\La}\dr_i^\La,
\qquad d_\La=d_{\la_r}\circ\cdots\circ d_{\la_1} \label{5.177}
\eeq
are total derivatives. There is the GDA epimorphism
$\cO^*_\infty [Y]\to \cO^*_\infty [U_Y]$ obtained as the restriction
of $\cO^*_\infty Y$ to the chart (\ref{jet1}).  Then $\cO^*_\infty
[Y]$ can be written in a coordinate form where the horizontal
one-forms $\{dx^\la\}$ and the contact one-forms
$\{\th^i_\La=dy^i_\La -y^i_{\la+\La}dx^\la\}$ are generating
elements of the $\cO^0_\infty [U_Y]$-algebra $\cO^*_\infty [U_Y]$.
Though $J^\infty Y$ is not a smooth manifold, the coordinate
transformations of elements of $\cO^*_\infty [Y]$ are smooth
since they are exterior forms on finite order jet
manifolds.

There is the canonical decomposition $\cO^*_\infty
[Y]=\oplus\cO^{k,m}_\infty [Y]$ of $\cO^*_\infty [Y]$ into
$\cO^0_\infty [Y]$-modules $\cO^{k,m}_\infty [Y]$ of $k$-contact
and $m$-horizontal forms together with the corresponding
projectors
\be
h_k:\cO^*_\infty [Y]\to \cO^{k,*}_\infty [Y], \qquad
h^m:\cO^*_\infty [Y]\to \cO^{*,m}_\infty [Y].
\ee
Accordingly, the exterior differential on $\cO_\infty^* [Y]$ is
split into the sum $d=d_H+d_V$ of the nilpotent total and vertical
differentials
\be
d_H(\f)= dx^\la\w d_\la(\f), \qquad d_V(\f)=\th^i_\La \w
\dr^\La_i\f, \qquad \f\in\cO^*_\infty [Y].
\ee

In particular, any finite order Lagrangian on a fiber bundle $Y\to
X$ is a density
\mar{0512}\beq
L=\cL\om \in \cO^{0,n}_\infty[Y], \qquad \om=dx^1\w\cdots\w dx^n,
\qquad n=\di X. \label{0512}
\eeq

In the framework of Lagrangian formalism, we deal with
differential operators of the following type. Let
\be
W\to Y\to X, \qquad Z\to Y\to X
\ee
be composite bundles, including $W=Y$, and let $Z\to Y$ be a
vector bundle. By a $k$-order differential operator on $W\to X$
taking its values in $Z\to X$ is throughout meant a bundle
morphism
\mar{k5}\beq
\Delta: J^kW\ar_Y Z. \label{k5}
\eeq
Its kernel $\Ker \Delta$ is defined as the inverse image of the
canonical zero section of $Z\to Y$. In an equivalent way, the
differential operator (\ref{k5}) is represented  by a section
$\Delta$ of the vector bundle $J^kW\op\times_Y Z\to J^kW$. Given
bundle coordinates $(x^\la,y^i,w^r)$ on $W$ and $(x^\la,y^i, z^A)$
on $Z$ with respect to the fiber basis $\{e_A\}$ for $Z\to Y$,
this section reads
\mar{k6}\beq
\Delta=\Delta^A(x^\la,y^i_\La, w_\La^B)e_A, \qquad  0\leq|\La|\leq
k. \label{k6}
\eeq
Then the differential operator (\ref{k5}) is also represented  by
an element
\mar{k7}\beq
\Delta=\Delta^A(x^\la,y^i_\La, w_\La^B)z_A \in \cO_\infty^0
[W\op\times Z^*] \label{k7}
\eeq
of the GDA $\cO_\infty^* [W\op\times_X Z^*]$, where $Z^*\to Y$ is
the dual of $Z\to Y$ with coordinates $(x^\la,y^i, z_A)$.

If $W\to Y$ is a vector bundle, a differential operator $\Delta$
(\ref{k5}) on the composite bundle $W\to Y\to X$ is said to be
linear if it is linear on the fibers of the vector bundle $J^kW\to
J^kY$. In this case, its representations (\ref{k6}) and (\ref{k7})
take the form
\mar{0608,'}\ben
&& \Delta=\op\sum_{0\leq|\Xi|\leq k} \Delta^{A,\Xi}_r(x^\la,
y^i_\La) w^r_\Xi e_A, \qquad 0\leq |\La|\leq k, \label{0608}\\
&& \Delta=\op\sum_{0\leq|\Xi|\leq k} \Delta^{A,\Xi}_r(x^\la,
y^i_\La) w^r_\Xi z_A, \qquad 0\leq |\La|\leq k. \label{0608'}
\een

In particular, every Lagrangian $L$ (\ref{0512}) defines the
Euler--Lagrange operator
\mar{k20}\beq
\dl L=\op\sum_{0\leq|\La|} (-1)^{|\La|}d_\La(\dr^\La_i
\cL)dy^i\w\om \label{k20}
\eeq
on $Y$ taking the values in the vector bundle
\mar{0548}\beq
V^*Y\op\ot_Y\op\w^n_X T^*X\to Y.
\label{0548}
\eeq
It is represented by the exterior form
\mar{0513}\beq
\dl L=\cE_i\th^i\w\om=\op\sum_{0\leq|\La|}
(-1)^{|\La|}d_\La(\dr^\La_i \cL)\th^i\w\om\in \cO^{1,n}_\infty[Y],
\label{0513}
\eeq
where
\mar{k8}\beq
\dl\f= \op\sum_{0\leq|\La|} (-1)^{\nm\La}\th^i\w
[d_\La(\dr^\La_i\rfloor d\f)], \qquad \f\in \cO^{*,n}_\infty [Y],
\label{k8}
\eeq
is the variational operator acting on $\cO^{*,n}_\infty [Y]$ so
that $\dl\circ d_H=0$ and $\dl\circ \dl=0$. There is the canonical
decomposition
\mar{+421}\beq
dL=\dl L-d_H\Xi, \label{+421}
\eeq
where $\Xi_L=L+\Xi$ is a Lepagean equivalent of $L$. It reads
\be
&& \Xi_L =L+\op\sum_{s=0}F^{\la\nu_s\ldots\nu_1}_i
\th^i_{\nu_s\ldots\nu_1}\w\om_\la, \\
&& F_i^{\nu_k\ldots\nu_1}= \dr_i^{\nu_k\ldots\nu_1}\cL-d_\la
F_i^{\la\nu_k\ldots\nu_1} +h_i^{\nu_k\ldots\nu_1}, \qquad
\om_\la=\dr_\la\rfloor\om,
\ee
where functions $h$ obey the relations $h^\nu_i=0$,
$h_i^{(\nu_k\nu_{k-1})\ldots\nu_1}=0$ \cite{got}.

\begin{rem} Given a Lagrangian $L$ and its Euler-Lagrange operator
$\dl L$ (\ref{0513}), we further abbreviate $A\ap 0$ with an
equality which holds on-shell. This means that $A$ is an element
of a module over the ideal $I_L$ of the ring $\cO^0_\infty [Y]$
which is locally generated by the variational derivatives $\cE_i$
and their total derivatives $d_\La\cE_i$. One says that $I_L$ is a
differential ideal because, if a local functions $f$ belongs to
$I_L$, then every total derivative $d_\La f$ does as well.
\end{rem}

\begin{rem}
We will use the relations
\mar{0606a-c}\ben
&& \op\sum_{0\leq |\La|\leq k}B^\La d_\La A'= \op\sum_{0\leq
|\La|\leq k} (-1)^{|\La|}d_\La (B^\La) A' +
d_H\si, \label{0606a}\\
&& \op\sum_{0\leq |\La|\leq k} (-1)^{|\La|}d_\La(B^\La A)=
\op\sum_{0\leq |\La|\leq k} \eta (B)^\La d_\La A, \label{0606b}
\\ && \eta (B)^\La = \op\sum_{0\leq|\Si|\leq
k-|\La|}(-1)^{|\Si+\La|} C^{|\Si|}_{|\Si+\La|} d_\Si B^{\Si+\La},
\qquad C^a_b=\frac{b!}{a!(b-a)!}, \label{0606c}
\een
for arbitrary exterior forms $A'\in \cO^{*,n}_\infty [Q]$,
$A\in\cO^*_\infty [Q]$ and local functions $B^\La\in \cO^0_\infty
[Q]$ on jet manifolds of a fiber bundle $Q\to X$. Since
$\op\sum_{a=0}^k(-1)^aC^a_k=0$ for $k>0$, it is easily verified
that
\mar{0606}\beq
(\eta\circ\eta)(B)^\La=B^\La. \label{0606}
\eeq
\end{rem}

\section{Gauge symmetries in a general setting}

Let $\gd\cO^0_\infty Y$ be the $\cO^0_\infty [Y]$-module of
derivations of the $\Bbb R$-algebra $\cO^0_\infty [Y]$. Any $\vt\in
\gd\cO^0_\infty Y$ yields a graded derivation (the interior
product) $\vt\rfloor\f$ of the GDA $\cO^*_\infty [Y]$ given by the
relations
\be
&&\vt\rfloor df=\vt(f), \qquad  f\in \cO^0_\infty [Y], \\
&& \vt\rfloor(\f\w\si)=(\vt\rfloor \f)\w\si
+(-1)^{|\f|}\f\w(\vt\rfloor\si), \qquad \f,\si\in \cO^*_\infty
[Y],
\ee
and a derivation $\bL_\vt$ (the Lie derivative) which satisfies the
conditions
\mar{0515,a,b}\ben
&& \bL_\vt\f=\vt\rfloor d\f+ d(\vt\rfloor\f), \qquad \f\in
\cO^*_\infty [Y], \label{0515}\\
&& \bL_\vt(\f\w\f')=\bL_\vt(\f)\w\f' +\f\w\bL_\vt(\f'),
\label{0515a}\\
&& \bL_\vt(d_H\f)=d_H(\bL_\vt\f). \label{0515b}
\een
Relative to an atlas (\ref{jet1}), a derivation
$\vt\in\gd\cO^0_\infty Y$ reads
\mar{g3}\beq
\vt=\vt^\la \dr_\la + \vt^i\dr_i + \op\sum_{|\La|>0}\vt^i_\La
\dr^\La_i, \label{g3}
\eeq
where the tuple of derivations $\{\dr_\la,\dr^\La_i\}$ is defined
as the dual of the set $\{dx^\la, dy^i_\La\}$ of generating
elements for the $\cO^0_\infty[Y]$-algebra $\cO^*_\infty[Y]$ with
respect to the interior product $\rfloor$ \cite{cmp}.

A derivation $\vt$ (\ref{g3}) is called contact if the Lie
derivative $\bL_\vt$ (\ref{0515}) preserves the contact ideal of
the GDA $\cO^*_\infty [Y]$ generated by contact forms. A
derivation $\vt$ (\ref{g3}) is contact iff
\mar{g4}\beq
\vt^i_\La=d_\La(\vt^i-y^i_\m\vt^\m)+y^i_{\m+\La}\vt^\m, \qquad
0<|\La|. \label{g4}
\eeq
Any contact derivation admits the canonical horizontal splitting
\mar{g5'}\beq
\vt=\vt_H +\vt_V=\vt^\la d_\la + (\up^i\dr_i + \op\sum_{0<|\La|}
d_\La \up^i\dr_i^\La), \qquad \up^i= \vt^i-y^i_\m\vt^\m.
\label{g5'}
\eeq
Its vertical part $\vt_V$  is completely determined by the first
summand
\mar{0641}\beq
\up=\up^i(x^\la,y^i_\La)\dr_i, \qquad 0\leq |\La|\leq k.
\label{0641}
\eeq
This is a section of the pull-back $VY\op\times_Y J^kY\to J^kY$ of
the vertical tangent bundle $VY\to Y$ onto $J^kY$ \cite{fat} and,
thus, it is a $k$-order $VY$-valued differential operator on $Y$.
One calls $\up$ (\ref{0641}) a generalized vector field on $Y$.

\begin{prop}  \label{g75} \mar{g75}
It follows from the splitting (\ref{+421}) that the Lie derivative
of a Lagrangian $L$ (\ref{0512}) along a contact derivation $\vt$
(\ref{g5}) fulfills the first variational formula
\mar{g8'}\beq
\bL_\vt L= \up\rfloor\dl L +d_H(h_0(\vt\rfloor\Xi_L)) +\cL d_V
(\vt_H\rfloor\om), \label{g8'}
\eeq
where $\Xi_L$ is a Lepagean equivalent of $L$ \cite{cmp}.
\end{prop}

A contact derivation $\vt$ (\ref{g5}) is called  variational if
the Lie derivative (\ref{g8'}) is $d_H$-exact, i.e., $\bL_\vt
L=d_H\si$, $\si\in \cO^{0,n-1}_\infty[Y]$. A glance at the
expression (\ref{g8'}) shows that: (i) a contact derivation $\vt$
is variational only if it is projectable onto $X$, (ii) $\vt$ is
variational iff its vertical part $\vt_V$ is well, (iii) it is
variational if $\up\rfloor \dl L$ is $d_H$-exact. By virtue of
item (ii), we can restrict our consideration to vertical contact
derivations
\mar{g5}\beq
\vt=\op\sum_{0\leq |\La|} d_\La \up^i\dr_i^\La. \label{g5}
\eeq
A generalized vector field $\up$ (\ref{0641}) is called a
variational symmetry of a Lagrangian $L$ if it generates a
variational vertical contact derivation (\ref{g5}).

A Lagrangian system on a fiber bundle $Y\to X$ is said to be a
gauge theory if its Lagrangian $L$ admits a family of variational
symmetries parameterized by elements of a vector bundle $E\to Y$
as follows.

Let $E\to Y$ be a vector bundle coordinated by
$(x^\la,y^i,\xi^r)$. Given a Lagrangian $L$ on $Y$, let us
consider its pull-back, say again $L$, onto $E$. Let $\vt_E$ be a
contact derivation of the $\Bbb R$-algebra $\cO^0_\infty [E]$, whose
restriction
\mar{0508}\beq
\vt=\vt_E|_{\cO^0_\infty Y}=
\op\sum_{0\leq|\La|}d_\La\up^i\dr_i^\La \label{0508}
\eeq
to $\cO^0_\infty [Y]\subset \cO^0_\infty [E]$ is linear in
coordinates $\xi^r_\Xi$. It is determined by a generalized vector
field $\up_E$ on $E$ whose canonical projection
\be
\up:J^kE\ar^{\up_E} VE\to E\op\times_Y VY
\ee
(see the exact sequence (\ref{0640}) below) is a linear
$VY$-valued differential operator
\mar{0509}\beq
\up= \op\sum_{0\leq|\Xi|\leq
m}\up^{i,\Xi}_r(x^\la,y^i_\Si)\xi^r_\Xi \dr_i \label{0509}
\eeq
on $E\to Y\to X$. Let $\vt_E$ be a variational symmetry of a
Lagrangian $L$ on $E$, i.e.,
\be
 \bL_{\vt_E}L=d_H\si.
\ee
Then one says that $\up$ (\ref{0509}) is a gauge symmetry of a
Lagrangian $L$.

\begin{rem}
Note that any generalized vector field $\up$ (\ref{0509}) gives
rise to a generalized vector field $\up_E$ on $E$ and, thus,
defines a contact derivation $\vt_E$ of $\cO^0_\infty [E]$.
Indeed, let us consider the exact sequence of vector bundles
\mar{0640}\beq
0\to V_YE\to VE\to E\op\times_Y VY\to 0, \label{0640}
\eeq
where $V_YE$ is the vertical tangent bundle of  $E\to Y$. Its
splitting $\G$ lifts $\up$ to the generalized vector field
$\up_E=\G\circ\up$ on $E$ such that the Lie derivative
\mar{k10}\beq
\bL_{\vt_E}L=\up\rfloor\dl L+ d_H(\vt\rfloor\Xi_L) \label{k10}
\eeq
depends only on $\up$, but not a lift $\G$.
\end{rem}

\begin{rem} If $\up$ (\ref{0509}) is a gauge symmetry, we obtain
from the first variational formula (\ref{k10}) the weak
conservation law
\mar{k11}\beq
0\ap d_H(\vt\rfloor\Xi_L-\si), \label{k11}
\eeq
where
\mar{k13}\beq
J=\vt\rfloor\Xi_L= \op\sum_{0\leq|\La|} J^{\la,\La}_r\xi^r_\La
\om_\la, \label{k13}
\eeq
is a Noether current.
\end{rem}

\section{Noether's second theorem}

Let us start with the notion of a Noether identity.

\begin{defi} \label{0568} \mar{0568}
Given a Lagrangian $L$ (\ref{0512}) and its Euler--Lagrange
operator $\dl L$ (\ref{0513}), let $E\to Y$ be a vector bundle and
$\Delta$ a linear differential operator of order $0\leq m$ on the
composite bundle (\ref{0548}) with the values in the density-dual
$\ol E^*$ (\ref{0630}) of $E$ which obeys the condition
(\ref{0549}). This condition is called the Noether identity, and
$\Delta$ is the Noether operator.
\end{defi}

Given bundle coordinates $(x^\la,y^i, \ol y_i)$ on the fiber
bundle (\ref{0548}) and  $(x^\la,y^i,\xi^r)$ on $E$, the Noether
operator $\Delta$ is represented by the
 density
\mar{0547}\beq
\Delta= \Delta_r\xi^r\om= \op\sum_{0\leq|\La|\leq m}
\Delta^{i,\La}_r(x^\la,y^j_\Si)\ol y_{\La i} \xi^r\om\in
\cO^{0,n}_\infty[E\op\times_Y V^*Y], \qquad 0\leq|\Si|\leq m.
\label{0547}
\eeq
Then the Noether identity (\ref{0549}) takes the coordinate form
\mar{0550}\beq
\op\sum_{0\leq|\La|\leq m} \Delta^{i,\La}_r d_\La \cE_i
\xi^r\om=0. \label{0550}
\eeq

\begin{theo} \label{0510} \mar{0510}
If a Lagrangian $L$ (\ref{0512}) admits a gauge symmetry $\up$
(\ref{0509}), its Euler--Lagrange operator obeys the Noether
identity (\ref{0550}) where the Noether operator (\ref{0547}) is
\mar{0511}\ben
&& \Delta=\eta(\up)=\op\sum_{0\leq|\Si|\leq m}
(-1)^{|\Si|}d_\Si(\up^{i,\Si}_r\ol y_i)\xi^r\om=
\op\sum_{0\leq|\La|\leq m} \eta(\up)^{i,\La}_r\ol y_{\La i}
\xi^r\om,
\label{0511}\\
&& \eta(\up)^{i,\La}_r=\op\sum_{0\leq|\Si|\leq
m-|\La|}(-1)^{|\Si+\La|}C^{|\Si|}_{|\Si+\La|} d_\Si
\up^{i,\Si+\La}_r. \nonumber
\een
Conversely, if the Euler--Lagrange operator of a Lagrangian $L$
obeys the Noether identity (\ref{0550}), this Lagrangian admits a
gauge symmetry  $\up$ (\ref{0509}) where
\mar{0638}\ben
&& \up=\eta(\Delta)=\op\sum_{0\leq|\Si|\leq m}
(-1)^{|\Si|}d_\Si(\Delta^{i,\Si}_r\xi^r)\dr_i=\op\sum_{0\leq|\La|\leq
m}\eta(\Delta)^{i,\La}_r\xi^r_\La \dr_i,
\label{0638}\\
&& \eta(\Delta)^{i,\La}_r=\op\sum_{0\leq|\Si|\leq
m-|\La|}(-1)^{|\Si+\La|} C^{|\Si|}_{|\Si+\La|}d_\Si
\Delta^{i,\Si+\La}_r.
\een
\end{theo}

\begin{proof} Given an operator $\up$ (\ref{0509}), the
operator $\Delta=\eta(\up)$ (\ref{0511}) is defined in accordance
with Theorem \ref{0610} in Appendix. Since the  density
\be
\up\rfloor\dl L=\up^i\cE_i\om=\op\sum_{0\leq|\Xi|\leq m}
\up^{i,\Xi}_r\xi^r_\Xi\cE_i\om
\ee
is $d_H$-exact, the Noether identity
\be
\dl(\up\rfloor\dl L)=\eta(\up)\circ \dl L=0
\ee
holds. Conversely, any operator $\Delta$ (\ref{0547}) defines the
generalized vector field $\up=\eta(\Delta)$ (\ref{0638}). Due to
the Noether identity (\ref{0550}), we obtain
\be
&&0=\op\sum_{0\leq|\La|\leq m}
\xi^r\Delta^{i,\La}_rd_\La\cE_i\om=\op\sum_{0\leq|\La|\leq m}
(-1)^{|\La|}d_\La(\xi^r \Delta^{i,\La}_r)\cE_i\om +
d_H\si=\\
&& \qquad \op\sum_{0\leq|\Xi|\leq m}\up^{i,\Xi}_r\xi^r_\Xi\cE_i\om
+d_H\si=\up\rfloor \dl L+d_H\si,
\ee
i.e., $\up$ is a gauge symmetry of $L$.
\end{proof}

By virtue of the relations (\ref{0634}), there is one-to-one
correspondence between gauge symmetries of a Lagrangian $L$ and
the Noether identities for $\dl L$.

\begin{ex} \label{0655} \mar{0655} If a gauge symmetry
\mar{0656}\beq
\up=(\up_r^i\xi^r +\up^{i,\m}_r\xi^r_\m)\dr_i \label{0656}
\eeq
is of first jet order in parameters,  the corresponding Noether
operator (\ref{0511}) and Noether identity take the form
\mar{0657,8}\ben
&& \Delta=[(\up^i_r -d_\m \up^{i,\m}_r)\ol y_i-
\up^{i,\m}_r\ol y_{\m i}]\xi^r\om,\label{0657}\\
&& [\up^i_r\cE_i - d_\m(\up^{i,\m}_r\cE_i)]\xi^r\om=0.
\label{0658}
\een
\end{ex}

Any Lagrangian $L$ has gauge symmetries. In particular, there
always exist trivial gauge symmetries
\be
\up=\op\sum_\La \eta(M)^{i,\La}_r\xi^r_\La, \qquad
M^{i,\La}_r=\op\sum_\Si T^{i,j,\La,\Si}d_\Si\cE_j,  \qquad
T_r^{j,i,\La,\Si}=-T_r^{i,j,\Si,\La},
\ee
corresponding to the trivial Noether identity
\be
\op\sum_{\Si,\La}T_r^{j,i,\La,\Si}d_\Si\cE_j d_\La\cE_i=0.
\ee
Furthermore, given a gauge symmetry $\up$ (\ref{0509}), let $E'\to
Y$ be a vector bundle and $h$ a linear differential operator on
some composite bundle $E'\to Y\to X$, coordinated by
$(x^\la,y^i,\xi'^s)$, with values in the vector bundle $E\to Y$.
Then the composition
\be
\up'=\up\circ h=\up'^{i,\La}_s\xi'^s_\La\dr_i, \qquad
\up'^{i,\La}_s=\op\sum_{\Xi+\Xi'=\La}\op\sum_{0\leq|\Si|\leq
m-|\Xi|} \up^{i,\Xi+\Si}_rd_\Si h^{r,\Xi'}_s,
\ee
is a variational symmetry of the pull-back of a Lagrangian $L$
onto $E'$, i.e., a gauge symmetry. In view of this ambiguity, we
agree to say that a gauge symmetry $\up$ (\ref{0509}) of a
Lagrangian $L$ is complete if any different gauge symmetry
$\up'_0$ of $L$ factorizes through $\up$ as
\be
\up'=\up\circ h + T, \qquad T\ap 0.
\ee
A complete gauge symmetry always exists, but the vector bundle of
its parameters need not be finite-dimensional.

Accordingly, given the Noether operator (\ref{0547}), let $H$ be a
linear differential operator on $\ol E^*\to Y\to X$ with values in
the density-dual $\ol E'^*$ (\ref{0630}) of some vector bundle
$E'\to Y$. Then the composition $\Delta'=H\circ\Delta$ is also a
Noether operator. We agree to call the Noether operator
(\ref{0547}) complete if a different Noether operator $\Delta'$
factors through $\Delta$ as
\be
\Delta'=H\circ\Delta + F, \qquad F\ap 0.
\ee

\begin{prop} \label{0646} \mar{0646}
A gauge symmetry $\up$ of a Lagrangian $L$ is complete iff 
the associated Noether operator is also.
\end{prop}

\begin{proof} The proof follows at once from Proposition
\ref{0621} in Appendix. Given a gauge symmetry $\up$ of $L$, let
$\up'$ be a different gauge symmetry. If $\eta(\up)$ is a complete
Noether operator, then
\be
\eta(\up')=H\circ\eta(\up) + F, \qquad F\ap 0,
\ee
and, by virtue of the relations (\ref{0622}), we have
\be
\up'=\up\circ\eta(H) +\eta(F),
\ee
where $\eta(F)\ap 0$ because $I_L$ is a differential ideal. The
converse is similarly proved.
\end{proof}

\section{Reducible gauge theories}

Let us extend Noether's second theorem to the analysis of
reducible gauge systems.

\begin{defi} \label{0573} \mar{0573} A complete Noether operator
$\Delta\not\ap 0$ (\ref{0547}) and the corresponding Noether
identity (\ref{0549}) are said to be {\it $N$-stage reducible}
($N=0,1,\ldots$) if there exist vector bundles $E_k\to Y$ and
differential operators $\Delta_k$, $k=0,\ldots,N$, such that:

(i)  $\Delta_k$ is a linear differential operator on the
density-dual $\ol E^*_{k-1}\to Y\to X$ of $E_{k-1}\to Y$ with
values in the density-dual $\ol E^*_k$ of $E_k$, where $E_{-1}=E$;

(ii) $\Delta_k\not\ap 0$ for all $k=0,\ldots,N$;

(iii) $\Delta_k\circ \Delta_{k-1}\ap 0$ for all $k=0,\ldots,N$,
where $\Delta_{-1}=\Delta$;

(iv) if $\Delta'_k$ is another differential operator possessing
these properties, then it factors through $\Delta_k$ on-shell.
\end{defi}

In particular, a zero-stage reducible Noether operator is called
reducible. In this case, given bundle coordinates
$(x^\la,y^i,\ol\xi_r)$ on $\ol E^*$ and $(x^\la,y^i,\xi^{r_0})$ on
$E_0$, a differential operator $\Delta_0$ reads
\mar{0576}\beq
\Delta_0=\op\sum_{0\leq|\Xi|\leq m_0} \Delta^{r,\Xi}_{r_0}
\ol\xi_{\Xi r}\xi^{r_0}\om. \label{0576}
\eeq
Then the reduction condition $\Delta_0\circ\Delta\ap 0$ takes the
coordinate form
\mar{0570}\beq
\op\sum_{0\leq|\Xi|\leq m_0}
\Delta^{r,\Xi}_{r_0}d_\Xi(\op\sum_{0\leq|\La|\leq
m}\Delta^{i,\La}_r\ol y_{\La i})\xi^{r_0}\om\ap 0, \label{0570}
\eeq
i.e., the left hand-side of this expression takes the form
\be
\op\sum_{0\leq|\Si|\leq m_0+m}M^{i,\Si}_{r_0} \ol y_{\Si
i}\xi^{r_0}\om,
\ee
where all the coefficients $M^{i,\Si}_{r_0}$ belong to the ideal
$I_L$.

\begin{defi} \label{0574} \mar{0574} A complete gauge
symmetry $\up\not\ap 0$ (\ref{0509}) is said to be $N$-stage
reducible if  there exist vector bundles $E_k\to Y$ and
differential operators $\up^k$, $k=0,\ldots,N$, such that:

(i)  $\up^k$ is a linear differential operator on the composite
bundle $E_k\to Y\to X$ with values in the vector bundle
$E_{k-1}\to Y$;

(ii) $\up^k\not\ap 0$ for all $k=0,\ldots,N$;

(iii) $\up^{k-1}\circ \up^k\ap 0$ for all $k=0,\ldots,N$, where
$\up^k$, $k=-1$, stands for $\up$;

(iv) if $\up'^k$ is another differential operator possessing these
properties, then $\up'^k$ factors through $\up^k$ on-shell.
\end{defi}

\begin{theo} \label{0580} \mar{0580}
A gauge symmetry $\up$ is $N$-stage reducible iff the
associated Noether identity is also.
\end{theo}

\begin{proof}
The proof follows at once from Theorem \ref{0610} and Proposition
\ref{0621} in Appendix. Let us put $\Delta_k=\eta(\up^k)$,
$k=0,\ldots,N$. If $\up^k\ap 0$, then $\eta(\up^k)\ap 0$ because
$I_L$ is a differential ideal. By the same reason, if $\up^{k-1}$
and $\up^k$ obey the reduction condition $\up^{k-1}\circ \up^k\ap
0$, then
\be
\eta(\up^{k-1}\circ \up^k)=\eta(\up^k)\circ\eta(\up^{k-1})\ap 0.
\ee
The converse  is justified in a similar way. The equivalence of
the conditions in items (iv) of Definitions \ref{0573} and
\ref{0574} is proved similarly to that in Proposition \ref{0646}.
\end{proof}

\begin{rem}
Let a gauge symmetry $\up$ (\ref{0509}) be reducible. Given bundle
coordinates $(x^\la, y^i, \xi^{r_0})$ on $E_0$, the differential
operator $\up^0$ reads
\be
\up^0=\op\sum_{0\leq|\La|\leq m_0}
\up^{r,\La}_{r_0}\xi^{r_0}_\La\dr_r,
\ee
and the reduction condition $\up\circ\up^0\ap 0$ takes the form
\be
\op\sum_{0\leq|\Xi|\leq m}\up^{i,\Xi}_r
d_\Xi(\op\sum_{0\leq|\La|\leq m_0} \up^{r,\La}_{r_0}\xi^{r_0}_\La)
\dr_i\ap 0.
\ee
In particular, it follows that the Noether current $J$ (\ref{k13})
vanishes on-shell if
\be
\xi^r= \op\sum_{0\leq|\La|\leq m_0} \up^{r,\Xi}_{r_0}\xi^{r_0}_\La
\ee
and, consequently, its components $J^{\la,\La}_r$ are not
independent.
\end{rem}

\section{Example I}

This example addresses the gauge model of principal connections on
a principal bundle $P\to X$ with a structure Lie group $G$ whose
automorphisms need not be vertical. In a general setting, the
gauge-natural prolongations of $P$ and the associated
natural-gauge bundles can be considered \cite{fat}.

Principal connections on a principal bundle $P\to X$ are
represented by sections of the quotient
\mar{0654}\beq
C=J^1P/G\to X,\label{0654}
\eeq
called the bundle of principal connections. This is an affine
bundle coordinated by $(x^\la, a^r_\la)$ such that, given a
section $A$ of $C\to X$, its components $A^r_\la=a^r_\la\circ A$
are coefficients of the familiar local connection form (i.e.,
gauge potentials). We consider the GDA $\cO^*_\infty [C]$.

Infinitesimal generators of one-parameter groups of automorphisms
of a principal bundle $P$ are $G$-invariant projectable vector
fields on $P\to X$. They are associated to sections of the vector
bundle $T_GP=TP/G\to X$. This bundle is endowed with the
coordinates $(x^\la,\tau^\la=\dot x^\la,\xi^r)$ with respect to
the fiber bases $\{\dr_\la, e_r\}$ for $T_GP$, where $\{e_r\}$ is
the basis for the right Lie algebra $\ccG$ of $G$ such that
$[e_p,e_q]=c^r_{pq}e_r.$ If
\mar{0652}\beq
 u=u^\la\dr_\la
+u^r e_r, \qquad v=v^\la\dr_\la +v^r e_r, \label{0652}
\eeq
are sections of $T_GP\to X$, their bracket reads
\be
[u,v]=(u^\m\dr_\m v^\la -v^\m\dr_\m u^\la)\dr_\la +(u^\la\dr_\la
v^r - v^\la\dr_\la u^r +c^r_{pq}u^pv^q)e_r.
\ee
Any section $u$ of the vector bundle $T_GP\to X$ yields the vector
field
\mar{0653}\beq
u_C=u^\la\dr_\la +(c^r_{pq}a^p_\la u^q +\dr_\la u^r -a^r_\m\dr_\la
u^\m)\dr^\la_r \label{0653}
\eeq
on the bundle of principal connections $C$ (\ref{0654})
\cite{book00}.

In order to describe a gauge symmetry in this gauge model, let us
consider the bundle product
\mar{0659}\beq
E=C\op\times_X T_GP, \label{0659}
\eeq
coordinated by $(x^\la, a^r_\la, \tau^\la, \xi^r)$. It can be
provided with the generalized vector field
\mar{0660}\beq
\up_E= \up=(c^r_{pq}a^p_\la \xi^q + \xi^r_\la
-a^r_\m\tau^\m_\la-\tau^\m a_{\m\la}^r)\dr^\la_r. \label{0660}
\eeq
With a subbundle $V_GP=VP/G\to X$ of the vector bundle $T_GP$
coordinated by $(x^\la, \xi^r)$, we have the exact sequence of
vector bundles
\be
0\to V_GP \ar T_GP\to TX \to 0.
\ee
The pull-back of this exact sequence via $C$ admits the canonical
splitting which takes the coordinate form
\mar{gr9}\beq
\tau^\la\dr_\la +\xi^r e_r = \tau^\la(\dr_\la +a^r_\la e_r) +
(\xi^r - \tau^\la a^r_\la)e_r. \label{gr9}
\eeq
Due to this splitting, the generalized vector field (\ref{0660})
is brought into the form
\mar{gr8}\beq
\up=(c^r_{pq}a^p_\la \xi'^q + \xi'^r_\la +
\tau^\m\cF_{\la\m}^r)\dr^\la_r, \qquad \xi'^r=\xi^r -\tau^\la
a^r_\la. \label{gr8}
\eeq

This generalized vector field is a gauge symmetry of the global
Chern--Simons Lagrangian in gauge theory on a principal bundle
with a structure semi-simple Lie group $G$ over a
three-dimensional base $X$. Given a section $B$ of $C\to X$ (i.e.,
a background gauge potential), this Lagrangian reads
\mar{r50}\ben
&& L= [\frac12a^G_{mn} \ve^{\al\bt\g}a^m_\al(\cF^n_{\bt\g}
-\frac13 c^n_{pq}a^p_\bt a^q_\g)  -\frac12a^G_{mn}
\ve^{\al\bt\g}B^m_\al(F(B)^n_{\bt\g} \label{r50}\\
&& \qquad -\frac13 c^n_{pq}B^p_\bt B^q_\g) -d_\al(a^G_{mn}
\ve^{\al\bt\g}a^m_\bt B^n_\g)]d^3x, \nonumber\\
&& F(B)^r_{\la\m}=\dr_\la B^r_\m-\dr_\m B^r_\la +c^r_{pq}B^p_\la
B^q_\m, \qquad \cF^r_{\la\m}=a^r_{\la\m}-a^r_{\m\la}
+c^r_{pq}a^p_\la a^q_\m, \nonumber
\een
where $a^G$ is the Killing form  \cite{bor,chern}. Its first term
is the well-known local Chern--Simons Lagrangian, the second one
is a density on $X$, and the Lie derivative of the third term is
$d_H$-exact due to the relation (\ref{0515b}). The corresponding
Noether identities (\ref{0658}) read
\mar{gr1,2}\ben
&& c^r_{pq}a^p_\la\cE_r^\la - d_\la(\cE_q^\la)=0,
\label{gr1}\\
&& -a^r_{\m\la}\cE^\la_r +d_\la(a^r_\m\cE^\la_r)=0. \label{gr2}
\een
The first one is the well-known Noether identity corresponding to
the vertical gauge symmetry
\be
 \up=(c^r_{pq}a^p_\la \xi^q +
\xi^r_\la)\dr^\la_r.
\ee
The second Noether identity (\ref{gr2}) is brought into the form
\be
-a^r_\m[c^r_{pq}a^p_\la\cE_r^\la - d_\la(\cE_q^\la)]
+\cF^r_{\la\m}\cE^\la_r=0,
\ee
i.e., it is equivalent to the Noether identity
$\cF^r_{\la\m}\cE^\la_r=0$, which also comes from the splitting
(\ref{gr9}) of the generalized vector field $\up$. This Noether
identity however is trivial since $\cF^r_{\la\m}=0$ is the kernel
of the Euler--Lagrange operator of the Chern--Simons Lagrangian
(\ref{r50}).

In order to obtain a gauge symmetry of the Yang--Mills Lagrangian,
one should complete the generalized vector field (\ref{0660}) with
the term acting on a world metric.

Let  $LX$ the fiber bundle of linear frames in the tangent bundle
$TX$ of $X$. It is a principal bundle with the structure group
$GL(n,\Bbb R)$, $n=\di X$, which admits reductions to its maximal compact
subgroup $O(n)$. Global sections of the quotient bundle
$\Si=LX/O(n)$ are Riemannian metrics on $X$. If $X$ obeys the
well-known topological conditions, pseudo-Riemannian metrics on
$X$ are similarly described. Being an open subbundle of the tensor
bundle $\op\vee^2 TX$, the fiber bundle $\Si$ is provided with
bundle coordinates $\si^{\m\nu}$. It admits the canonical lift
\mar{gr3}\beq
u_\Si=u^\la\dr_\la +(\si^{\nu\bt}\dr_\nu u^\al
+\si^{\al\nu}\dr_\nu u^\bt)\frac{\dr}{\dr \si^{\al\bt}}
\label{gr3}
\eeq
of any vector field $u=u^\la\dr_\la$ on $X$. We describe the gauge
system of principal connections and a dynamic metric field on the
bundle product
\mar{gr10}\beq
E=C\op\times_X \Si\op\times_X T_GP, \label{gr10}
\eeq
coordinated by $(x^\la,a^r_\la,\si^{\al\bt}, \tau^\la, \xi^r)$. It
can be provided with the generalized vector field
\mar{gr11}\beq
\up=(c^r_{pq}a^p_\la \xi^q + \xi^r_\la -a^r_\m\tau^\m_\la-\tau^\m
a_{\m\la}^r)\dr^\la_r + (\si^{\nu\bt}\tau_\nu^\al +\si^{\al\nu}
\tau_\nu^\bt-\tau^\la\si_\la^{\al\bt})\frac{\dr}{\dr
\si^{\al\bt}}. \label{gr11}
\eeq
This is a gauge symmetry of the sum of the Yang--Mills Lagrangian
$L_{\rm YM}(\cF^r_{\al\bt}, \si^{\m\nu})$ and a Lagrangian of a
metric field. The corresponding Noether identities read
\be
&& c^r_{pq}a^p_\la\cE_r^\la - d_\la(\cE_q^\la)=0,\\
&& -a^r_{\m\la}\cE^\la_r +d_\la(a^r_\m\cE^\la_r)
-\si_\m^{\al\bt}\cE_{\al\bt} - 2d_\nu (\si^{\nu\bt}
\cE_{\mu\bt})=0.
\ee
The first one is the Noether identity (\ref{gr1}). Then the second
identity is brought into the form
\be
\cF^r_{\la\m}\cE^\la_r - 2\nabla_\nu (\si^{\nu\bt}
\cE_{\mu\bt})=0,
\ee
where $\nabla_\nu$ are covariant derivatives with respect to the
Levi--Civita connection
\be
K=dx^\la\ot (\dr_\la + K_\la{}^\m{}_\nu \dot x^\nu\dot\dr_\mu),
\qquad K_\la{}^\m{}_\nu =-\frac12\si^{\nu\bt}(\si_{\la\bt\m} +
\si_{\m\bt\la} -\si_{\bt\la\m}).
\ee

\section{Example II}

Let us consider gauge theory of skew symmetric tensor fields.
These are exterior forms on a base manifold $X$ of degree more
than one. We need not specify a gauge model, but refer to the
topological BF theory \cite{birm}. This is a theory of two
exterior forms $A$ and $B$ of form degree $|A|=\di X-|B|-1$.
Another example is a gauge theory of an exterior form $A$ in the
presence of a background metric on $X$ whose Lagrangian is similar
to that of an electromagnetic field.

A generic gauge system of skew symmetric tensor fields is defined
on the fiber bundle
\mar{k30}\beq
Y=\op\w^p_XT^*X\op\oplus_X \op\w^q_XT^*X, \label{k30}
\eeq
coordinated by $(x^\la, A_{\m_1\ldots\m_p},B_{\nu_1\ldots\nu_q})$.
The corresponding GDA is $\cO^*_\infty[Y]$. There are the
canonical $p$- and $q$-forms
\mar{k31}\beq
A=\frac{1}{p!}A_{\m_1\ldots\m_p}dx^{\m_1}\w\cdots\w dx^{\m_p}\in
\cO^{0,p}_\infty[Y], \quad
B=\frac{1}{q!}B_{\nu_1\ldots\nu_q}dx^{\nu_1}\w\cdots\w dx^{\nu_q}
\in \cO^{0,q}_\infty[Y] \label{k31}
\eeq
on $Y$. A Lagrangian of the above-mentioned topological BF theory
reads
\mar{k32}\beq
L_{\rm BF}=A\w d_HB, \qquad p+q=n-1. \label{k32}
\eeq

A gauge symmetry of a generic gauge system of skew symmetric
tensor fields, e.g., of the Lagrangian (\ref{k32}) is the
following. Let us consider the fiber bundle
\mar{k33}\beq
E=Y\op\times_X \op\w^{p-1}_XT^*X\op\times_X \op\w^{q-1}_XT^*X,
\label{k33}
\eeq
where
\be
\op\w^{p-1}T^*X\op\times_X \op\w^{q-1}T^*X
\ee
is the fiber bundle of gauge parameters with coordinates $(x^\la,
\ve_{\m_1\ldots\m_{p-1}},\xi_{\nu_1\ldots\nu_{q-1}})$. Let
\be
&& \ve=\frac{1}{(p-1)!}\ve_{\m_1\ldots\m_{p-1}}dx^{\m_1}\w\cdots\w
dx^{\m_{p-1}}\in \cO^{0,p-1}_\infty[E], \\
&& \xi=\frac{1}{(q-1)!}\xi_{\nu_1\ldots\nu_{q-1}}dx^{\nu_1}\w\cdots\w
dx^{\nu_{q-1}} \in \cO^{0,q-1}_\infty[E]
\ee
be canonical exterior forms like (\ref{k31}) on $E$ (\ref{k33}).
The above mentioned gauge symmetry is given by the generalized
vector field
\mar{k35}\beq
\up=d_{\m_1}\ve_{\m_2\ldots\m_p}\frac{\dr}{\dr A_{\m_1\ldots\m_p}}
+ d_{\nu_1}\xi_{\nu_2\ldots\nu_q}\frac{\dr}{\dr
B_{\nu_1\ldots\nu_q}}, \label{k35}
\eeq
which acts on the exterior forms $A$ and $B$ (\ref{k31}) by the
law
\mar{k36}\beq
\bL_{\vt_E}A=d_H\ve, \qquad \bL_{\vt_E}B=d_H\xi. \label{k36}
\eeq
In accordance with the formula (\ref{0658}), the corresponding
Noether identity takes the form
\mar{k37}\beq
-\ve_{\m_2\ldots\m_p}d_{\m_1} \cE^{\m_1\ldots\m_p}=0, \qquad
-\xi_{\nu_2\ldots\nu_q}d_{\nu_1} \cE^{\nu_1\ldots\nu_q}=0.
\label{k37}
\eeq

For instance, the equalities
\be
&& \bL_{\vt_E}L_{\rm BF}= (\bL_{\vt_E}A)\w d_HB) +A\w
(\bL_{\vt_E}d_HB)= d_H\ve\w d_HB + A\w (\bL_{\vt_E}d_HB)\\
&& \qquad = d_H(\ve\w d_HB)
\ee
show that the generalized vector field (\ref{k35}) is a gauge
symmetry of the Lagrangian $L_{\rm BF}$ (\ref{k32}). This
Lagrangian provides the Euler--Lagrange equations
\be
d_HA=0, \qquad d_HB=0,
\ee
and the Noether identity (\ref{k37}) is brought into the form
\be
d_Hd_HA\equiv 0, \qquad d_Hd_H B\equiv 0.
\ee
It should be emphasized that the gauge symmetry (\ref{k35}) by no
means exhausts all variational symmetries of the Lagrangian
$L_{\rm BF}$. Any generalized vector field $\up$ on $Y$ such that
$\bL_\vt A$ is $d_H$-exact and $\bL_\vt B$ is $d_H$-closed is a
variational symmetry of this Lagrangian.

The gauge symmetry $\up$ (\ref{k35}) is reducible. Without a loss
of generality, let us put $q\geq p$. Then $\up$ is $q$-stage
reducible as follows. Let us consider vector bundles
\mar{k40}\ben
&& E_k=Y\op\times_X \op\w^{p-k-2}_XT^*X\op\times_X
\op\w^{q-k-2}_XT^*X, \qquad 0\leq k< p-2, \nonumber\\
&& E_k=Y\op\times_X \Bbb R \op\times_X
\op\w^{q-p}_XT^*X, \qquad k=p-2, \label{k40}\\
&& E_k=Y\op\times_X \op\w^{q-k-2}_XT^*X, \qquad k>p-2, \nonumber
\een
over $Y$ provided with fiber coordinates
\be
(\ve^k_{\m_1\ldots\m_{p-k-2}}, \xi^k_{\nu_1\ldots\nu_{q-k-2}}),
\qquad (\al,\xi^{p-2}_{\nu_1\ldots\nu_{q-p}}), \qquad
(\xi^k_{\nu_1\ldots\nu_{q-k-2}}),
\ee
respectively. Then the differential operators
\be
&& \up_0=d_{\m_1}\ve^0_{\m_2\ldots\m_{p-1}}\frac{\dr}{\dr
\ve_{\m_1\ldots\m_{p-1}}} +
d_{\nu_1}\xi^0_{\nu_2\ldots\nu_{q-1}}\frac{\dr}{\dr
\xi_{\nu_1\ldots\nu_{q-1}}}, \\
&& \up_k=d_{\m_1}\ve^k_{\m_2\ldots\m_{p-k-1}}\frac{\dr}{\dr
\ve^{k-1}_{\m_1\ldots\m_{p-k-1}}} +
d_{\nu_1}\xi^k_{\nu_2\ldots\nu_{q-k-1}}\frac{\dr}{\dr
\xi^{k-1}_{\nu_1\ldots\nu_{q-k-1}}}, \qquad 0< k< p-2,\\
&& \up_{p-2}=d_\m\al \frac{\dr}{\dr \ve^{p-3}_\m} +
d_{\nu_1}\xi^{p-2}_{\nu_2\ldots\nu_{q-p+1}}\frac{\dr}{\dr
\xi^{p-3}_{\nu_1\ldots\nu_{q-p+1}}},\\
&& \up_k= d_{\nu_1}\xi^k_{\nu_2\ldots\nu_{q-k-1}}\frac{\dr}{\dr
\xi^{k-1}_{\nu_1\ldots\nu_{q-k-1}}}, \qquad k>p-2,
\ee
satisfy the conditions of Definition \ref{0574}.

Accordingly, there is a family of associated $k$-stage Noether
operators $\Delta_k$. Let the density duals $\ol E^*\to Y$, $\ol
E^*_k\to Y$ of the vector bundles $E\to Y$ (\ref{k33}), $E_k\to Y$
(\ref{k40}) be provided with fiber coordinates
\be
(\ol\ve^{\m_1\ldots\m_{p-1}}, \ol\xi^{\nu_1\ldots\nu_{q-1}}),
\qquad (\ol\ve_k^{\m_1\ldots\m_{p-k-2}},
\ol\xi_k^{\nu_1\ldots\nu_{q-k-2}}), \qquad
(\ol\al,\ol\xi_{p-2}^{\nu_1\ldots\nu_{q-p}}), \qquad
(\ol\xi_k^{\nu_1\ldots\nu_{q-k-2}}),
\ee
respectively. Then we have
\be
&& \Delta_0 =-[\ve^0_{\m_2\ldots\m_{p-1}}d_{\m_1}
\ol\ve^{\m_1\ldots\m_{p-1}} +
\xi^0_{\nu_2\ldots\nu_{q-1}}d_{\nu_1}
\ol\xi_{k-1}^{\nu_1\ldots\nu_{q-k-1}}]\om,\\
&& \Delta_k =-[\ve^k_{\m_2\ldots\m_{p-k-1}}d_{\m_1}
\ol\ve_{k-1}^{\m_1\ldots\m_{p-k-1}} +
\xi^k_{\nu_2\ldots\nu_{q-k-1}}d_{\nu_1}
\ol\xi_{k-1}^{\nu_1\ldots\nu_{q-k-1}}]\om, \qquad 0< k< p-2,\\
&& \Delta_{p-2} =-[\al d_\m \ol\ve_{p-3}^\m +
\xi^{p-2}_{\nu_2\ldots\nu_{q-p+1}}d_{\nu_1}
\ol\xi_{p-3}^{\nu_1\ldots\nu_{q-p+1}}]\om,\\
&& \Delta_k = - \xi^k_{\nu_2\ldots\nu_{q-k-1}}d_{\nu_1}
\ol\xi_{k-1}^{\nu_1\ldots\nu_{q-k-1}}\om, \qquad k>p-2.
\ee

\section{Appendix. Differential operators on dual fiber bundles}

Given a fiber bundle $Y\to X$, let $E\to Y$, $Q\to Y$ be vector
bundles coordinated by $(x^\la,y^i,\xi^r)$ and $(x^\la,y^i,q^a)$,
respectively. Let $E^*$, $Q^*$ be their duals and  $\ol E^*$, $\ol
Q^*$ their density-duals (\ref{0630}) coordinated by
$(x^\la,y^i,\xi_r)$, $(x^\la,y^i,q_a)$ and  $(x^\la,y^i,\ol
\xi_r)$, $(x^\la,y^i,\ol q_a)$, respectively. Let $\up$ be a
linear $Q$-valued differential operator on $E\to Y\to X$. It is
represented by the function (\ref{0608'}):
\mar{0614}\beq
\up=\up^aq_a=\op\sum_{0\leq|\La|\leq m}
\up^{a,\La}_r(x^\la,y^i_\Si) \xi^r_\La q_a\in
\cO^0_\infty[E\op\times_Y Q^*], \qquad 0\leq|\Si|\leq m.
\label{0614}
\eeq
Let $\Delta$ be an $m$-order linear differential operator on the
density-dual $\ol Q^*\to Y\to X$ of $Q\to Y$
 with values in the density-dual
$\ol E^*$ of $E$. It is represented by the function
\mar{k1}\beq
\ol\Delta=\op\sum_{0\leq|\La|\leq m}
\Delta^{a,\La}_r(x^\la,y^i_\Si)\ol q_{\La a}\ol\xi^r \in
\cO^0_\infty[(\ol E^*)^*\op\times_Y \ol Q^*], \qquad
0\leq|\Si|\leq m. \label{k1}
\eeq
Let $J\om$ be a volume form on $X$ such that $\ol\xi_r=J\xi_r$ and
$\xi^r=J\ol\xi^r$. Then the function (\ref{k1}) defines the
density $\Delta=\ol\Delta J\om$ which reads
\mar{0615}\beq
 \Delta=\Delta_r\xi^r\om=\op\sum_{0\leq|\La|\leq m}
\Delta^{a,\La}_r(x^\la,y^i_\Si)\ol q_{\La a}\xi^r\om\in
\cO^{0,n}_\infty[E\op\times_Y \ol Q^*], \qquad 0\leq|\Si|\leq m.
\label{0615}
\eeq

\begin{theo} \label{0610} \mar{0610}
Any linear $Q$-valued differential operator $\up$ (\ref{0614}) on
$E\to Y\to X$ yields the linear $\ol E^*$-valued differential
operator
\mar{0616}\ben
&& \eta(\up)= \op\sum_{0\leq|\Si|\leq m}
(-1)^{|\Si|}d_\Si(\up^{a,\Si}_r\ol
q_a)\xi^r\om=\op\sum_{0\leq|\La|\leq m}
\eta(\up)^{a,\La}_r \ol q_{\La a}\xi^r\om, \label{0616}\\
&&\eta(\up)^{a,\La}_r=\op\sum_{0\leq|\Si|\leq
m-|\La|}(-1)^{|\Si+\La|}C^{|\Si|}_{|\Si+\La|} d_\Si
(\up^{a,\Si+\La}_r), \nonumber
\een
on $\ol Q^*\to Y\to X$.  Conversely, any linear $\ol E^*$-valued
differential operator $\Delta$ (\ref{0615}) on $\ol Q^*\to Y\to X$
defines the linear $Q$-valued  differential operator
\mar{0617}\ben
&&\eta(\Delta)= \op\sum_{0\leq|\Si|\leq m}
(-1)^{|\Si|}d_\Si(\Delta^{a,\Si}_r\xi^r)q_a=\op\sum_{0\leq|\La|\leq
m} \eta(\Delta)^{a,\La}_r
\xi^r_\La q_a, \label{0617}\\
&&\eta(\Delta)^{a,\La}_r=\op\sum_{0\leq|\Si|\leq
m-|\La|}(-1)^{|\Si+\La|}C^{|\Si|}_{|\Si+\La|} d_\Si
(\Delta^{a,\Si+\La}_r), \nonumber
\een
on $E\to Y\to X$. The relations (\ref{0634}) hold.
\end{theo}

\begin{proof} One must show that the differential operators given
by the local coordinate expressions (\ref{0616}) and (\ref{0617})
are globally defined. The function $\up$ (\ref{0614}) yields the
density
\mar{0637}\beq
\ol\up=\op\sum_{0\leq|\La|\leq m} \up^{a,\La}_r \xi^r_\La \ol
q_a\om\in \cO^{0,n}_\infty[E\op\times_Y Q^*].\label{0637}
\eeq
Its Euler--Lagrange operator
\be
\dl(\ol\up)=\cE_i dy^i\w\om +\cE_r d\xi^r\w\om + \cE^ad\ol
q_a\w\om
\ee
takes its values in the fiber bundle
\mar{0631}\beq
V^*(E\op\times_Y Q^*)\op\ot_{E\op\times_Y Q^*}\op\w^n_X T^*X,
\label{0631}
\eeq
where $V^*(E\op\times_Y Q^*)$ is the vertical cotangent bundle of
the fiber bundle $E\op\times_Y Q^*\to X$. There is its canonical
projection
\mar{0633}\beq
\al_E: V^*(E\op\times_Y Q^*)\to V^*E\to V^*_YE, \label{0633}
\eeq
onto the vertical cotangent bundle $V^*_YE$ of $E\to Y$. Then we
obtain a differential operator $(\al_E\circ\dl)(\ol\up)$ on
$E\op\times_Y Q^*$ with values in the fiber bundle
$V^*_YE\op\ot_E\op\w^n T^*X$. It reads
\be
(\al_E\circ\dl)(\ol\up)=\cE_r\ol
d\xi^r\ot\om=\op\sum_{0\leq|\La|\leq m}
(-1)^{|\La|}d_\La(\up^{a,\La}_r \ol q_a)\ol d\xi^r\ot\om,
\ee
where $\{\ol d\xi^r\}$ is the fiber basis for $V^*_YE\to E$ and the
tensor product $\ot$ is over $C^\infty(X)$. Due to the canonical
isomorphism
$V^*_YE=E^*\op\times_Y E$, this operator defines the  density
(\ref{0616}). Conversely,  the Euler--Lagrange operator of the density
(\ref{0615}) takes its values in the fiber bundle (\ref{0631}) and reads
\mar{0632}\beq
\dl(\Delta)=\cE_i dy^i\w\om +\cE_r d\xi^r\w\om + \cE^a d\ol
q_a\w\om. \label{0632}
\eeq
In  order to repeat the above mentioned procedure, let us consider
a volume form $J\om$ on $X$ and substitute $d\ol
q_a\w\om=Jdq_a\w\om$ into the expression (\ref{0632}). Using the
projection
\be
\al_Q: V^*(E\op\times_Y Q^*)\to V^*_YQ^*
\ee
like $\al_E$ (\ref{0633}) and the canonical isomorphism
$V^*_YQ^*=Q\op\times_Y Q^*$, we come to the density
\be
\op\sum_{0\leq|\La|\leq m} (-1)^{|\La|} d_\La(\Delta^{a,\La}_r
\xi^r) q_a J\om \in\cO^{0,n}_\infty[E\op\times_Y Q^*]
\ee
and, hence, the function (\ref{0617}). The relations (\ref{0634})
result from the relation (\ref{0606}).
\end{proof}

Relations (\ref{0634}) show that the intertwining operator $\eta$
(\ref{0616}) -- (\ref{0617}) provides a  bijection between the
sets Diff$(E,Q)$ and Diff$(\ol Q^*,\ol E^*)$ of differential
operators (\ref{0614}) and (\ref{0615}).

\begin{prop} \label{0621} \mar{0621}
Compositions of operators $\up\circ\up'$ and $\Delta'\circ \Delta$
obey the relations (\ref{0622}).
\end{prop}

\begin{proof}
It suffices to prove the first relation. Let $\up\circ\up'\in{\rm
Diff}(E',Q)$ be a composition of differential operators
$\up\in{\rm Diff}(E,Q)$ and $\up'\in{\rm Diff}(E',E)$. Given
fibred coordinates $(\xi^r)$ on $E\to Y$, $(\e^p)$ on $E'\to Y$
and $(\ol q_a)$ on $\ol Q^*\to Y$, this composition defines the
density (\ref{0637})
\be
\ol{\up\circ\up'}=\op\sum_\La \up^{a,\La}_r
d_\La(\op\sum_\Si\up'^{r,\Si}_p\e^p_\Si)\ol q_a\om.
\ee
Following the relation (\ref{0606a}), one can bring this density
into the form
\be
\op\sum_\Si\up'^{r,\Si}_p\e^p_\Si\op\sum_\La (-1)^{|\La|} d_\La
(\up^{a,\La}_r \ol q_a)\om + d_H\si=
\op\sum_\Si\up'^{r,\Si}_p\e^p_\Si\op\sum_\La\eta(\up)^{a,\La}_r
\ol q_{\La a}\om +d_H\si.
\ee
Its Euler--Lagrange operator projected to $V^*_YE'\op\ot_{E'}
T^*X$ is
\be
\op\sum_\Si (-1)^{|\Si|} d_\Si(\up'^{r,\Si}_p
\op\sum_\La\eta(\up)^{a,\La}_r \ol q_{\La a})\ol
d\e^p\ot\om=\op\sum_\Si\eta(\up')^{r,\Si}_p
d_\Si(\op\sum_\La\eta(\up)^{a,\La}_r \ol q_{\La a})\ol
d\e^p\ot\om,
\ee
that leads to the desired composition $\eta(\up')\circ\eta(\up)$.
\end{proof}


\begin{thebibliography}{ddd}



\bibitem{ander} Anderson I 1992 Introduction to the variational
bicomplex {\it Contemp. Math.} {\bf 132}  51-73


\bibitem{barn} Barnich G, Brandt F and Henneaux M 2000 Local
BRST cohomology in gauge theories {\it Phys. Rep.} {\bf 338}
439-569

\bibitem{bash} Bashkirov D and Sardanashvily G 2005 On the BV
quantization of gauge gravitation theory {\it Int. J. Geom.
Methods Mod. Phys.} {\bf 2} No.2

\bibitem{birm} Birmingham D, Blau M, Rakowski M and Thompson G
1991 Topological field theories {\it Phys. Rep.} {\bf 209} 129-340

\bibitem{bor} Borowiec A, Ferraris M and Francaviglia M 2003 A
covariant formalism of Chern--Simons gravity {\it J. Phys. A} {\bf
36} 2589-2598

\bibitem{fat} Fatibene L and Francaviglia M 2003 {\it Natural and
Gauge Natural Formalism for Classical Field Theories. A Geometric
Perspective Including Spinors and Gauge Theories} (Dordrecht:
Kluwer)


\bibitem{fish} Fish M and Henneaux M 1990 Homological
perturbation theory and algebraic structure of the
antifield-antibracket formalism for gauge theories {\it Commun.
Math. Phys.} {\bf 128} 627-640

\bibitem{fulp} Fulp R, Lada T and Stasheff J 2003 Noether
variational Theorem II and the BV formalism {\it Rend. Circ. Mat.
Palermo (2) Suppl.} No. 71 115-126

\bibitem{book00} Giachetta G, Mangiarotti L and
Sardanashvily G 2000 {\it Connections in Classical and Quantum
Field Theory} (Singapore: World Scientific)

\bibitem{jmp} Giachetta G, Mangiarotti L and
Sardanashvily G 2001  Cohomology of the infinite-order jet space
and the inverse problem {\it J. Math. Phys.} {\bf 42} 4272-4282

\bibitem{chern} Giachetta G, Mangiarotti L and Sardanashvily G 2003
Noether conservation laws in higher-dimensional Chern-Simons theory {\it
Mod. Phys. Lett. A} {\bf 18} 2645-2651

\bibitem{cmp} Giachetta G, Mangiarotti L and
Sardanashvily G 2005 Lagrangian supersymmetries depending on
derivatives. Global analysis and cohomology. {\it Commun. Math.
Phys.} (accepted) ({\it Preprint} hep-th/0407185)

\bibitem{gom} Gomis J, Par\'\i s J and Samuel S 1995
Antibracket, antifields and gauge theory quantization {\it Phys.
Rep} {\bf 295} 1-145

\bibitem{got} Gotay M 1991 A multisymplectic framework for
classical field theory and the calculus of variations. In: {\it
Mechanics, Analysis and Geometry: 200 Years after Lagrange}
(Amsterdam: North Holland) 203-235.

\bibitem{tak2} Takens F 1979 A global version of the inverse
problem of the calculus of variations {\it J. Diff. Geom.} {\bf
14} 543-562


\end{thebibliography}
\end{document}